\documentclass[12pt, a4paper]{article}
\newtheorem{thm}{{\bf T}{\footnotesize \bf HEOREM}}
\newtheorem{lm}[thm]{{\bf L}{\footnotesize \bf EMMA}}

\newtheorem{pro}[thm]{{\bf P}{\footnotesize \bf ROPOSITION}}
\newtheorem{conj}[thm]{{\bf C}{\footnotesize \bf ONJECTURE}}

\newenvironment{Proof}{\noindent {\it Proof.} }{\hbox{\rule{6pt}{6pt}}
\bigskip}
\newenvironment{Proofof}[1]
{\noindent {\it Proof of {#1}.} }{\hbox{\rule{6pt}{6pt}} \bigskip}

\usepackage[utf8]{inputenc}
\usepackage[dvipdfmx]{graphicx}
\usepackage{tikz}
\usepackage{tikzinclude}
\usepackage{amsmath, amsfonts}
\usepackage{makecell}

\begin{document}
\sloppy
\title{Connectivity and matching extendability of optimal $1$-embedded graphs on the torus}
\author{
Shohei Koizumi\thanks{%
Graduate School of Science and Technology,
Niigata University,
8050 Ikarashi 2-no-cho, Nishi-ku, Niigata, 950-2181, Japan.
E-mail: {\tt s-koizumi@m.sc.niigata-u.ac.jp}} \
and Yusuke Suzuki\thanks{%
Department of Mathematics, Niigata University, 
8050 Ikarashi 2-no-cho, Nishi-ku, Niigata, 950-2181, Japan.
Email: {\tt y-suzuki@math.sc.niigata-u.ac.jp}}
}
\date{}
\maketitle

\begin{abstract}
\noindent
In this paper, we discuss 
optimal $1$-toroidal graphs (abbreviated as O1TG), 
which are drawn on the torus so that 
every edge crosses another edge at most once, and has $n$ vertices and exactly $4n$ edges. 
We first consider connectivity of O1TGs, and 
give the characterization of O1TGs having 
connectivity exactly $k$ 
for each $k\in \{4, 5, 6, 8\}$. 
In our argument, we also show that there exists no O1TG 
having connectivity exactly $7$. 
Furthermore, using the result above, 
we discuss extendability of matchings, and 
give the characterization of $1$-, $2$- and $3$-extendable 
O1TGs in turn. 
\end{abstract}


\section{Introduction}\label{sect:intro}
\noindent
All graphs considered in this paper are finite, simple and connected. 
We denote the vertex set and the edge set of a graph $G$ 
by $V(G)$ and $E(G)$, respectively. 
The {\em order\/} of $G$ means the number of vertices of $G$. 
A path $P$ of length $k$ is a {\em $k$-path\/}, which 
has $k$ edges.  
In particular, if $k=0$, then $P$ is a trivial graph. 
A cycle of length $k$ is a {\em $k$-cycle\/}. 
For a cycle $C$ in a graph $G$, an edge $e \in E(G)$ such that 
$V(e) \subset V(C)$ and $e \notin E(C)$ is called a {\em chord\/} of $C$. 
A cycle $C$ in $G$ is {\em separating\/} 
if $G-V(C)$ is a disconnected graph. 
We denote the induced subgraph of $G$ by $S \subset V(G)$ by $G[S]$. 

A graph $G$ is {\em $1$-embeddable\/} on a closed surface $F^2$
if it can be drawn on $F^2$ 
so that every edge of $G$ crosses another edge at most once. 
The drawn image of $G$ on $F^2$ is a $1$-{\em embedded} graph on $F^2$. 
(We implicitly consider {\em good drawings\/}, that is, (i) vertices are on different 
points on the surface, (ii) no adjacent edges cross, (iii) no three edges cross 
at the same point, and (iv) any non-adjacent edges do not touch tangently.)
The study of $1$-{\em planar\/} graphs, which are 
$1$-embeddable graphs on the plane or the sphere, 
was first introduced by Ringel \cite{Ringel},
and recently developed in various points of view (see e.g., \cite{KLM, Book}); 
the drawn image is called a $1$-{\em plane graph\/}. 
It is known that if $G$ is a $1$-embedded graph on $F^2$ 
with at least three vertices,
then $|E(G)| \leq 4|V(G)| - 4\chi(F^2)$ holds, 
where $\chi(F^2)$ stands for the Euler characteristic of $F^2$ 
(see \cite{Nagasawa} for example).
In particular, a $1$-embedded graph $G$ on $F^2$ that 
satisfies the equality, that is $|E(G)| = 4|V(G)| - 4\chi(F^2)$, 
is {\em optimal\/}. 
An edge in a $1$-embedded graph $G$ is {\em crossing\/} 
if it crosses another edge, 
and {\em non-crossing\/} otherwise. 

It was shown in \cite{Nagasawa} that 
every simple optimal $1$-embedded graph $G$ on $F^2$ is obtained from a polyhedral 
quadrangulation $H$ by adding a pair of crossing edges in each face of $H$. 
We call the quadrangulation $H$, which consists of all the non-crossing edges of $G$, 
the {\em quadrangular subgraph\/} of $G$, 
and denote it by $Q(G) (=H)$. 
By the property above, every vertex of $G$ has even degree, that is, 
$G$ is Eulerian. 
For example, ``optimal $1$-plane graph'' is abbreviated as ``O1PG'' in past research; see e.g., 
\cite{O1PG, O1PM}.  
Similarly, ``optimal $1$-embedded graph on the torus (resp., the projective plane)'' 
is also said like ``{\em optimal $1$-toroidal graph\/}'' 
(resp., ``{\em optimal $1$-projective plane graph\/}'', 
and abbreviated as ``O1TG'' (resp., ``O1PPG'').

In the first half of the paper, we discuss connectivity of O1TGs. 
In \cite{FSS}, Fujisawa et. al discussed connectivity of O1PGs, 
and show that every O1PG has connectivity either $4$ or $6$; 
they actually gave the characterization of O1PGs having connectivity $k$ for each $k\in \{4,6\}$. 
Furthermore in \cite{Ko}, connectivity of O1PPGs were discussed, and 
it was shown that there exist O1PPGs that have connectivity exactly $5$.  
That is, every O1PPG has connectivity $4$, $5$ or $6$, 
and the authors 
actually characterized such O1PPGs with connectivity $k$ for 
each $k\in \{4, 5, 6\}$. 
Since every O1PG or O1PPG $G$ has a vertex of degree exactly $6$, 
and hence $\kappa(G)\leq 6$, where $\kappa(G)$ represents the connectivity of $G$.  
However for O1TGs, the situation is different, since there are infinitely many 
$8$-regular O1TGs whose quadrangular subgraphs are $4$-regular. 
The first main result in this paper is the characterization of $8$-connected O1TGs 
as follows, where $Q(p,q,r)$ is a well-known representation of $4$-regular quadrangulations 
in topological graph theory, 
which will be introduced in Section~\ref{sect:Q(G)}. 

\begin{thm}\label{thm:8-conn}
Let $G$ be an $8$-regular O1TG. 
Then, we have $\kappa(G) \in \{6, 8\}$.
Furthermore, $G$ is $8$-connected if and only if $Q(G)$ is not isomorphic to 
$Q(p, r, 3)$ $(p \geq 4$ and $r \geq 0)$. 
\end{thm}
%
%

Actually, it is easy to show that every optimal $1$-embedded graph on a closed 
surface is $4$-connected; this fact is proved as Theorem~\ref{thm:4-conn-gene} in this paper.  
Hence, we characterize O1TGs, which are not $8$-regular, 
having connectivity exactly $k \in \{4, 5\}$ as follows. 
(As a result, O1TGs other than those have connectivity exactly $6$.)


\begin{thm}\label{thm:4,5-conn}
Let $G$ be an O1TG. Then the followings hold: 
\begin{enumerate}
\item[(i)] $\kappa(G)=4$ if and only if $Q(G)$ has a trivial $4$-cycle that 
does not bound a face of $Q(G)$.   
\item[(ii)] Assume that $G$ is $5$-connected. Then $\kappa(G)=5$ if and only if $Q(G)$ has two homotopic cycles $xy_1z_1x$ and $xy_2z_2x$, 
where $x, y_1, y_2, z_1$ and $z_2$ are distinct vertices. 
\end{enumerate}
\end{thm}

In the latter half of the paper, we discuss matching extendability of 
O1TGs using the result related to connectivity of those graphs discussed above. 
A matching $M\subset E(G)$ of $G$ is {\em extendable\/} 
if $G$ has a perfect matching containing $M$. 
Moreover, 
a graph $G$ with at least $2m+2$ vertices is {\em $m$-extendable\/} 
if any matching $M$ in $G$ with $|M| = m$ is extendable. 
Matching extendability has been widely studied in literature (e.g., see \cite{P08}).
In particular, matching extendability of graphs embedded on closed surfaces was investigated in 
\cite{AKP, AP11, KNPS,P1988}; for example, 
it was proven as a basic result that no planar graph is $3$-extendable.

The matching extendability of $1$-embedded graphs on $F^2$ was first addressed  
in \cite{FSS}, and the authors proved that 
every O1PG of even order is $1$-extendable. 
Further in the same paper, they discussed $2$-extendability of O1PGs, and 
proved that an O1PG $G$ of even order is $2$-extendable
unless $G$ contains a barrier $4$-cycle, 
where a {\em barrier $k$-cycle\/} is a $k$-cycle of $Q(G)$ that bounds a $2$-cell containing 
odd number of vertices; it is indeed defined in Section~\ref{sect:Proofs}. 
As mentioned above, every O1PG has a vertex $v$ of degree exactly $6$, 
and it is clear that every O1PG is not $3$-extendable; consider three independent 
edges covering $N_G(v)$. 
In the same paper, they actually characterize three independent edges 
that are extendable. 
Moreover in \cite{Ko}, the argument was extended to O1PPGs. 

We first give the characterization of $1$-extendable O1TGs as follows, 
where $T^2$ represents the torus throughout the paper. 

\begin{thm}\label{thm:1-ext}
Let $G$ be an O1TG of even order.  
Then $G$ is $1$-extendable if and only if 
$Q(G)$ does not contain a subgraph $H$ satisfying the following conditions: 
\begin{enumerate}
\item[(i)]  $H$ is a quadrangulation of $T^2$, and 
\item[(ii)] every facial cycle of $H$ corresponds to a barrier $4$-cycle of $G$. 
\end{enumerate} 
\end{thm}

The following statement concern with the $2$-extendability of O1TGs 
is similar to those for O1PG's and O1PPGs. 

\begin{thm}\label{thm:2-ext}
An O1TG $G$ of even order 
is $2$-extendable if and only if  
$G$ has no barrier $4$-cycle. 
\end{thm}

Unlike in the case of O1PGs and O1PPGs, there exist $3$-extendable O1TGs, 
and those graphs are characterized as follows.   

\begin{thm}\label{thm:3-ext}
An O1TG $G$ of even order 
is $3$-extendable if and only if  
$G$ is $8$-regular. 
\end{thm}

This paper is organized as follows.
In the next section, 
we first define terminology used in the paper, and 
introduce the fundamental results holding for optimal $1$-embedded 
graphs on general closed surfaces. 
Next, we 
discuss connectivity of O1TGs 
and separating short cycles in quadrangular subgraphs.  
In Section~\ref{sect:Q(G)}, 
we provide the characterization of $8$-regular O1TGs 
having connectivity exactly $8$. 
Furthermore in Section~\ref{sect:not8}, 
we characterize not $8$-regular O1TGs 
having connectivity $k$ for each $k\in \{4, 5, 6\}$. 
In Section~\ref{sect:Proofs}, 
we discuss matching 
extendability of O1TGs, and show the characterization of 
$m$-extendable O1TGs for each $m \in \{1,2,3\}$.


\section{Preliminaries and basic results}\label{sect:pre}

A vertex set $S$ of a connected graph $G$ is a {\em cut\/} if 
$G - S$ has at least two connected components. 
A cut $S$ of $G$ is {\em minimal\/} if any proper subset of $S$ 
is not a cut of $G$. 
For a cut $S$ of $G$, if $|S| = k$, then we call $S$ a {\em $k$-cut\/} 
of $G$. 
We denote the number of odd components of 
$G-S$ for $S\subset V(G)$ by $C_o(G-S)$; that is, a connected component of 
odd order.

Let $G$ be a graph embedded on a closed surface $F^2$. 
Then a connected component of $F^2 - G$, which is as a topological space, 
is a {\em face\/} of $G$, 
and we denote the face set of $G$ by $F(G)$;  
that is, ``a face'' in this paper is not necessarily 
homeomorphic to an (open) $2$-cell. 
In particular, 
if every face of $G$ is homeomorphic to a $2$-cell, then 
$G$ is a {\em $2$-cell embedding\/} or {\em $2$-cell embedded graph\/} on $F^2$. 

In general, each boundary component of a face $f$ forms 
a closed walk of $G$. 
That is, the boundary of $f$, which is denoted by 
$\partial f$, is the union of closed walks of $G$. 
In particular, if $f$ has the unique boundary component, 
then it is said to be the {\em boundary closed walk\/} of $f$. 
Let $f$ be a face of $G$ embedded on $F^2$, and 
assume that $\partial f = W_1\cup \cdots \cup W_l$, 
where $W_i$ is a closed walk corresponding to a 
boundary component of $f$ for each $i\in \{1, \ldots, l\}$. 
Then, the invariant $\deg(f)$, which is called the {\em size\/} of $f$, 
is defined as follows where $|W_i|$ is the length of a closed walk $W_i$: 
$\deg(f) = \sum_{i=1}^{l}|W_i|$. 
Furthremore, we put $V(\partial f) = V(W_1) \cup \cdots \cup V(W_l)$. 
A $k$-{\em face\/} $f$ of $G$ is a face with 
$\deg(f) = k$. 
In particular, if a $k$-face $f$ is bounded by a closed walk $W_1=v_0v_1\cdots v_{k-1}v_0$ of length $k$, 
that is, $\partial f=W_1$, then 
$f$ is a $k$-{\em gonal face\/} of $G$. 
In this case, we simply denote $f=v_0v_1\cdots v_{k-1}$ in our latter arguments. 
Furthremore, if $f$ is a $2$-cell face, we say that $f$ is a $k$-{\em gonal\/} $2$-{\em cell face\/}. 
Observe that every $4$-face is a $4$-gonal face since it cannot have at least two boundary 
components, except for a very special case where the boundary consists of two independent 
edges such that the $4$-gonal face is incident to those edges on both sides. 

A simple closed curve $\gamma$ on a closed surface $F^2$ is 
{\em trivial\/} if $\gamma$ bounds a $2$-cell on $F^2$, and 
{\em essential\/} otherwise. 
We apply these definition to cycles of graphs embedded on $F^2$, 
regarding them as simple closed curves. 
A simple closed curve $\gamma$ on a closed surface $F^2$ 
is {\em surface separating\/} if $F^2 - \gamma$ is disconnected as a topological space. 
We also apply the definition to cycles of graphs on $F^2$. 
It is well-known that every surface separating simple closed curve 
on the sphere, the projective plane or the torus is trivial. 
The following proposition is known in 
topological graph theory, and is commonly used.

\begin{pro}[\cite{Naka}]
\label{prop:parity}
Let $G$ be a graph $2$-cell 
embedded on a closed surface $F^2$ so that each face is bounded 
by a closed walk of even length. Then the length of two cycles in $G$ have the same parity 
if they are homotopic to each other on $F^2$. 
Furthermore, there is no surface separating odd cycle in $G$. 
\end{pro}

%

We often consider the induced subgraph of 
$Q(G)$ by a cut $S$ of $G$ in our argument, 
which is $Q(G)[S]$ under our definition. 
However, when the underlying graph $G$ is clear, we use $Q[S]$ in place of $Q(G)[S]$, to simplify the notation. 
In the following five lemmas, we assume that 
$G$ is an optimal $1$-embedded graph on $F^2$, 
and $S\subset V(G)$ is a cut of $G$. 
First of all, we show the following one, which is related to Proposition~\ref{prop:parity}. 

\begin{lm}\label{lm:parity}
For every face $f$ of $Q[S]$, the size of $f$ is even. 
\end{lm}

\begin{Proof}
Consider the graph $H= Q(G) \cap (f\cup \partial f)$ embedded on the 
surface $F_0$ with the boundary components $W_1 \cup \cdots \cup W_l$. 
We glue a disc $D_i$ to $W_i$ so that the boundary of $D_i$ coincides with 
$W_i$ for each $i \in \{1, \ldots, l\}$, to obtain a closed surface $\tilde{F}_0$. 
Now, take the dual $H^*$ of $H$ on $\tilde{F}_0$. 
If the size of $f$ is odd, then it immediately contradicts the odd point theorem. 
\end{Proof}

The remaining four lemmas are actually proved in $\cite{Ko}$.

\begin{lm}[\cite{Ko}]\label{lm:sep}
Every face of $Q[S]$ contains at most 
one connected component of $G - S$. 
\end{lm}


\begin{lm}[\cite{Ko}]\label{lm:degree}
If $S$ is minimal, then the minimum degree of $Q[S]$ is at least $2$.
\end{lm}

\begin{lm}[\cite{Ko}]\label{lm:Q[S]}
If $Q[S]$ has $p$ faces each of which has the size 
at least $2q \geq 6$, then the following inequalities hold: 
\begin{enumerate}
\item[(i)] $|E(Q[S])| \geq 2|F(Q[S])| + (q -2)p$ 
\item[(ii)] $|S| - \chi(F^2) + (2-q)p \geq |F(Q[S])|$ 
\end{enumerate} 
\end{lm}

\begin{lm}[\cite{Ko}]\label{lm:edge-bound}
If $|S| \leq C_o(G-S)+2m$ holds for some integer $m$, 
then we have the following: 
$$2|F(Q[S])| + 2m - \chi(F^2) \geq |E(Q[S])|$$ 
\end{lm}

\section{Connectivity of $8$-regular O1TGs}\label{sect:Q(G)}

In this section, we discuss the connectivity of $8$-regular O1TGs. 
First of all, we show that every O1TG is $6$-connected. 

\begin{lm}\label{lm:6-conn.}
Let $G$ be an $8$-regular O1TG. 
Then $6 \leq \kappa(G) \leq 8$ holds. 
\end{lm}

\begin{Proof}
An $8$-regular O1TG $G$ has a $6$-regular triangulation $T$ as a spanning subgraph, 
which is obtained from $Q(G)$ by adding a diagonal edge in each face; 
note that $Q(G)$ is a $4$-regular quadrangulation of $T^2$. 
Then $T$ is $6$-connected by the result in \cite{Negami}, $\kappa(G) \geq 6$ holds. 
In addition, we clearly have $\kappa(G) \leq 8$, 
since every vertex of $G$ has degree exactly $8$. 
\end{Proof}


The following lemma describes the inner structures 
of $2$-cell regions of quadrangulations of closed surfaces bounded by 
closed walks of 
length $4$, $6$ and $8$. 
Note that a {\em region\/} in this paper might contain 
vertices and edges of the graph.   

\begin{lm}[\cite{SuzukiPG}]\label{lm:bcycle}
Let $G$ be a quadrangulation of a closed surface $F^2$ and let $D$ be 
a $2$-cell region bounded by a closed walk $C$ of length $4$, $6$ or $8$ 
such that 
\begin{enumerate}
\item[(i)] there is at least one vertex inside D,  
\item[(ii)] all vertices inside $D$ have degree at least $3$ and 
\item[(iii)] $D$ does not have a unique vertex $x$ of degree $4$ 
such that $N_G(x) \subset V(C)$ and $D$ contains exactly four faces incident to $x$ 
(when $|C| = 8$). 
\end{enumerate}

\noindent
Then, there exists a vertex of degree $3$ inside $D$.  
\end{lm}


The following lemma holds for quadrangulations with minimum degree $4$.

\begin{lm}\label{lm:4regu.}
Let $G$ be a quadrangulation of $F^2$ with minimum degree at least $4$ and let 
$D$ be a $2$-cell region bounded by a closed walk $C$ of length $4$, $6$ or $8$. 
Then the following holds: 
\begin{enumerate}
\item[(i)] If $|C| = 4$ or $6$, then there are no vertices inside $D$. 
\item[(ii)] If $|C|=8$, then there is the unique vertex of degree $4$ inside $D$. 
\end{enumerate}
\end{lm}

\begin{Proof}
It immediately follows from Lemma~\ref{lm:bcycle}. 
\end{Proof}

Next, we discuss $Q[S]$ for a $k$-cut $S$ in $8$-regular O1TGs, and 
present the following three lemmas, which are keys to prove our main result. 

%
%

\begin{lm}\label{lm:6cut2}
Let $G$ be an $8$-regular O1TG, and let $S$ be a minimal $k$-cut of $G$ with $k\in \{6, 7\}$. 
Then $Q[S]$ has at least two faces that are not homeomorphic to $2$-cells. 
\end{lm}

\begin{Proof}
Suppose, for a contradiction, that 
$Q[S]$ has at most one face that is not homeomorphic to $2$-cells. 
Since $S$ is a cut of $G$, 
$Q[S]$ has a $2$-cell face $f_1$, that contains a 
connected component of $G-S$ by Lemma~\ref{lm:sep}. 
Furthermore, $\deg(f_1)$ is even by Lemma~\ref{lm:parity}. 
Then, we have $\deg(f_1) \geq 8$ by (i) of Lemma~\ref{lm:4regu.}. 
If $\deg(f_1)=8$, then $f_1$ contains the 
unique vertex $v$ of $G$ by (ii) of Lemma \ref{lm:4regu.}. 
Since $k\in \{6, 7\}$, $\partial f$ is not a cycle.  
This contradicts that $G$ is simple since $\deg_G(v)=8$. 
Therefore, we have $\deg(f_1)\geq 10$. 

Next, we discuss the other face $f_2$ that contains a connected 
component of $G$; $f_2$ might not be a $2$-cell face. 
By the minimality of $S$, $V(\partial f_2)=S$. 
(Observe that $f_2$ contains at least one vertex of $G$ since 
no edge of $Q(G)$ other than those on the boundary is inside $f_2$;  
recall that $Q[S]$ is the induced subgraph of $Q(G)$.) 
Hence $f_2$ is a $k'$-face where $k'$ is an even number with at least $k$ 
by Lemma~\ref{lm:parity}. 

Then the inequality 
$2|E(Q[S])| \geq 4(|F(Q[S])|-2)+k' +10 = 4|F(Q[S])|+k' +2$ holds. 
On the other hand, we have $k-|E(Q[S])|+|F(Q[S])|\geq 0$, 
and hence $k+|F(Q[S])|\geq |E(Q[S])|$ holds. 
By combining the former inequality, we obtain 
$2k-(k' +2)\geq 2|F(Q[S])|$. 
Since $|F(Q[S])| \geq 2$, the equality in the inequality above must hold, 
and it implies that 
$Q[S]$ is a $2$-cell embedding. 
Thus, $Q[S]$ has exactly two faces $f_1$ and $f_2$ 
with $\deg(f_1)=10$ and $\deg(f_2)\in \{6,8\}$. 
However, $\deg(f_2)$ must be at least $10$ by the argument above, a contradiction. 
\end{Proof}


\begin{lm}\label{lm:6cut3}
Let $G$ be an $8$-regular O1TG, and let $S$ be a $6$-cut of $G$. 
Then $Q[S]$ is the union of two essential $3$-cycles 
that are homotopic to each other. 
\end{lm}

\begin{Proof}
By Lemma~\ref{lm:6cut2}, $Q[S]$ has at least two faces $f_1$ and $f_2$ that 
are not homeomorphic to $2$-cells. 
Assume that $\partial f_1 = W_1\cup \cdots \cup W_l$, where 
each $W_i$ is a boundary component that is a closed walk in $Q[S]$. 
First, we suppose $l=1$. 
In this case, $\partial f_1$ consists of the unique boundary component, 
and $f_1$ contains a handle. 
Under the condition, $f_2$ that is not a $2$-cell face cannot exist. 
Thus, we conclude that $l\geq 2$. 
 
Now we take a simple closed curve $\gamma_i$ along $W_i$ inside $f_1$ for 
each $i \in \{1, \ldots, l\}$. 
Suppose that $\gamma_i$ is trivial on $T^2$ for some $i\in \{1, \ldots, l\}$, 
say $\gamma_1$ without loss of generality.  
In this case, note that $V(W_1) \cap V(W_i)=\emptyset$ for each $i \neq 1$. 
If the interior of $\gamma_1$ contains a vertex other than those of $W_1$, 
then, $W_1$ would become a smaller cut set, contrary to $S$ being minimal. 
On the other hand, if the interior of $\gamma_1$ does not contain any vertex 
other than those of $W_1$, then it also concludes a contradiction since 
$V(W_2)\cup \cdots \cup V(W_l)$ is a smaller cut set of $G$. 
As a result, $\gamma_i$ is not trivial on $T^2$ for each $i\in \{1, \ldots, l\}$. 

Since $\gamma_1$ and $\gamma_2$ do not have any intersection, 
they are homotopic to each other on $T^2$. 
That is, $\gamma_1$ can be shifted through $f_1$ and aligned to $\gamma_2$, 
and the trajectory itself forms an annular region. 
Furthermore, this implies that no other boundary component cannot exist in 
the annular region, and hence we have $l=2$; 
observe that if there exists $W_3$ in the annular region, then $\gamma_3$ 
must be trivial, a contradiction. 
That is, the annular region is actually $f_1$. 

Since $G$ is simple, we have $V(W_1) \geq 3$ and $V(W_2) \geq 3$. 
Furthermore, if $V(W_1) \cap V(W_2) \neq \emptyset$, 
then the other non-2-cell face $f_2$ does not exist on $T^2$; 
observe that the same argument also holds for $f_2$, and 
we can take an essential simple closed curve in $f_2$. 
Thus, we have $|V(W_1)| = |V(W_2)| = 3$. 
Furthermore, each of $W_1$ and $W_2$ is a $3$-cycle under the situation, 
and hence we have got our desired conclusion. 
Observe that $f_2$ also has the boundary components 
$W_1$ and $W_2$, which are $3$-cycles. 
\end{Proof}


\begin{lm}\label{lm:min7-cut}
Every $8$-regular O1TG has no minimal $7$-cut. 
\end{lm}

\begin{figure}[t] 
\begin{center}
\input{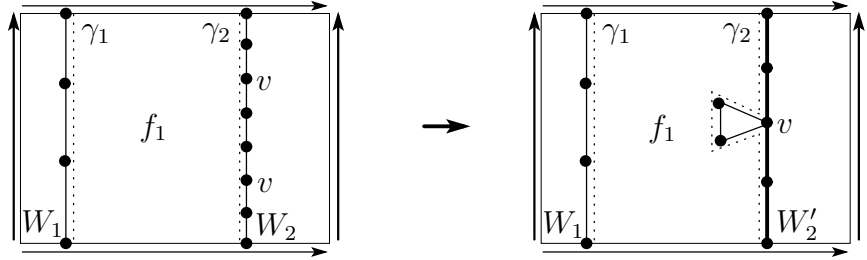}
\caption{Homotopic closed  walks $W_1, W_2$ and $W_2'$.}
\label{fig:homotopic} 
\end{center} 
\end{figure}

\begin{Proof}
Suppose, for a contradiction, that an $8$-regular O1TG $G$ has a 
minimal $7$-cut $S$. 
Then, by Lemma~\ref{lm:6cut2}, $Q[S]$ has at least two faces $f_1$ and $f_2$ that 
are not homeomorphic to $2$-cells. 
Actually, the same argument as that in the previous proof holds, 
and hence, we have $\partial f_1 = W_1 \cup W_2$ with $V(W_1) \geq 3$ and $V(W_2) \geq 3$, 
and $V(W_1) \cap V(W_2) = \emptyset$. 
Without loss of generality, we may assume that $|V(W_1)|=3$, that is, $W_1$ is a cycle of length $3$. 
If $W_2$ is not a cycle, then $W_2$ contains a closed walk $W_2'$ homotopic to $W_2$ 
with $|V(W_2)|>|V(W_2')|$ (see Figure~\ref{fig:homotopic}.)  
In this case, $V(W_1) \cup V(W_2')$ would become a smaller cut, contrary $S$ being minimal. 
Thus, $W_2$ is also a cycle. 
By Proposition~\ref{prop:parity}, we have $|W_2|=3$, and hence it contradicts $|S|=7$.    
\end{Proof}


As follows, $4$-regular quadrangulations of $T^2$ is completely classified in \cite{AL}, 
where $Q(p, r, q)$ is a commonly used notation that 
represents those graphs.  
(Every $4$-regular quadrangulations of $T^2$ can be cut along edges 
into a $p\times q$ rectangle, and the parameter $r$ represents 
the deviation when reconstructing the torus from the annulus with 
two boundary cycles of length $q$. 
See \cite{AL, Naka, NakaNega, Negami, Negami0} in details.)
In particular, $q$ represents the length of the meridian cycle of 
$Q(p, r, q)$. 
Note that to keep $Q(p, r, q)$ simple, we need the condition $q \geq 3$ at least. 

\begin{thm}[\cite{AL}]\label{thm:4reguq}
Every $4$-regular quadrangulation of $T^2$ is isomorphic to $Q(p, r, q)$ for 
some integeres $p \geq 1$, $q \geq 3$ and $r \geq0$. 
\end{thm}

%
%


Now we prove our first main result in the paper. 

\bigskip

\begin{Proofof}{Theorem \ref{thm:8-conn}}
By Lemmas~\ref{lm:6-conn.} and \ref{lm:min7-cut}, 
we have $\kappa(G)\in \{6,8\}$. 
Then, we show the latter part of the theorem below. 
The necessity is trivial, and hence we discuss the sufficiency of the statement. 
Let $G$ be an $8$-regular O1TG that is not $8$-connected. 
As above, we discuss the case of $\kappa(G)=6$. 
Then there exists a $6$-cut $S$ of $G$. 
By Lemma~\ref{lm:6cut3}, 
$Q[S]$ is the union of two essential $3$-cycles 
that are homotopic to each other, 
and hence $G$ is isomorphic to 
$Q(p, r, 3)$ $(p \geq 4$ and $r \geq 0)$; 
observe that each annular region of $T^2-Q[S]$ contains 
vertices of $G$, otherwise $S$ would not become a cut of $G$.  
\end{Proofof}


\section{Connectivity of not $8$-regular O1TGs}\label{sect:not8}

In this section, we discuss the connectivity of not $8$-regular O1TGs. 
At the beginning of the section, we prove the following theorem. 
(For the cases of the sphere and the projective plane, 
that is, for O1PGs and O1PPGs, the proofs are provided separately 
in \cite{FSS} and \cite{Ko}, respectively.)

\begin{thm}\label{thm:4-conn-gene}
Every optimal $1$-embedded graph on a closed surface is $4$-connected. 
\end{thm}

\begin{Proof}
Let $G$ be an optimal $1$-embedded graph on a closed surface $F^2$. 
Suppose, for a contradiction, that $G$ has a $3$-cut $S$. 
Then, we have $|F(Q[S])| \geq 2$ by Lemma~\ref{lm:sep}. 
Further, we have $|E(Q[S])| \geq 2|F(Q[S])|$ by Lemma~\ref{lm:parity}. 
This is a contradiction since we clearly have $|E(Q[S])|\leq 3$. 
\end{Proof}


The following is also a generalization of the results 
in \cite{FSS, Ko}. 

\begin{lm}\label{lm:4cycle}
Let $G$ be an optimal $1$-embedded graph on $F^2$, and let 
$S$ be a $4$-cut of $G$. 
Then, one of the followings holds: 
\begin{enumerate}
\item[(I)] $Q[S]$ has exactly three $4$-gonal faces 
each of which contains a vertex of $G$.  
In this case, $Q[S]$ has exactly six edges, 
and $F^2$ is nonorientable. 
\item[(II)] $Q[S]$ has exactly two faces each of which contains 
a vertex of $G$. Furthermore, one of them is a $4$-gonal face. 
\end{enumerate}
\end{lm}

\begin{Proof}
We put $S=\{s_1, s_2, s_3, s_4\}$. 
By Lemma~\ref{lm:sep}, we have $|F(Q[S])| \geq 2$. 
Furthermore, $|E(Q[S])| \geq 2|F(Q[S])|$ by Lemma~\ref{lm:parity}, 
and hence we obtain $|F(Q[S])| \leq 3$; recall that $|V(Q[S])|=|S|=4$, 
and then $|E(Q[S])|\leq 6$.  
First, we consider the case of $|F(Q[S])|=3$. 
In this case, we immediately have $|E(Q[S])|=6$ by the inequality above, 
and this implies that $Q[S]$ has exactly three $4$-faces; 
that is, three $4$-gonal faces. 
Let $f_1, f_2$ and $f_3$ denote such three $4$-faces, and 
we may assume that $f_1=s_1s_2s_3s_4$ without 
loss of generality. 
We may further assume that $s_1s_2$ is shared by two faces 
$f_1$ and $f_2$. 
If $f_2=s_1s_2s_3s_4$, then $Q[S]$ must be a $4$-cycle on the sphere 
by the simplicity of $G$, contradicting the existence of $f_3$. 
Thus, we have $f_2=s_1s_2s_4s_3$, and $F^2$ contains a M\"obius band 
that crosses edges $s_1s_2$ and $s_3s_4$.  
This is actually (I) in the statement. 

Next, we consider the case of $F(Q[S])=2$. 
Let $f_1$ and $f_2$ be two faces of $Q[S]$. 
Suppose, for a contradiction, that 
$\deg(f_1)=\deg(f_2)=6$ with $|E(Q[S])|=6$; the other cases 
clearly admit a $4$-face by a similar argument 
as above. 
Now, we consider the following two cases: 
Case (i) $\partial f_1$ has at least two boundary components, 
and Case (ii) $\partial f_1$ has exactly one boundary component. 

We first 
discuss Case (i). 
In this case, we may assume that $\partial f_1=W_1\cup W_2$ with 
$V(W_1)=\{s_1, s_2, s_3\}$ and $V(W_1)=\{s_1, s_2, s_4\}$. 
Under the situation, a $4$-cycle $s_1s_3s_2s_4s_1$ cuts off 
the region obtained from $f_1$ by identifying $s_1s_2$ from 
the other part, which is nothing but $f_2$, contradicting that $f_2$ 
has size $6$. 

Next, we assume Case (ii). 
By Lemma~\ref{lm:degree} and the simplicity of $G$, 
each $s_i$ appear on $\partial f_1$ at most twice. 
Then, we may assume that $\partial f_1=s_1s_2s_3s_1s_2s_4s_1$ 
without loss of generality; observe that the element of 
$S$ that appears on $\partial f_1$ twice must be on the 
diagonal position. 
Similar to the previous case, 
the cycle $s_2s_3s_1s_4s_2$ is a surface separating $4$-cycle;
that is, $f_2=s_2s_3s_1s_4$ is a $4$-gonal face of $Q[S]$, a contradiction.  
\end{Proof}


Next, we discuss minimal $5$-cuts in O1TGs. 

\begin{lm}\label{lm:5conn}
Let $G$ be an O1TG, and let 
$S$ be a minimal $5$-cut of $G$. 
Then $Q[S]$ has a $6$-face $f$ such that 
$V(\partial f)=S$ and $f$ contains a connected component of $G-S$. 
\end{lm}

\begin{Proof}
By Lemma~\ref{lm:parity}, 
every face of $Q[S]$ has even size. 
Note that since $S$ is minimal, 
$Q[S]$ has no $4$-face that contains a connected component of $G-S$. 
By Lemma~\ref{lm:sep}, $|F(Q[S])| \geq 2$ holds, and hence 
$Q[S]$ has at least two faces with size at least $6$. 
Then the inequality 
$|F(Q[S])| \leq 3$ 
follows from (ii) of Lemma~\ref{lm:Q[S]} with $|S|=5$, $\chi(F^2)=0$, $p\geq 2$ and 
$q\geq 3$. 
By combining $|F(Q[S])|\leq 3$ and Euler's formula $5-|E(Q[S])|+|F(Q[S])| \geq 0$, 
we obtain $|E(Q[S])|\leq 8$. 

First, we consider the case of $|E(Q[S])|=8$. 
Then all the equalities of inequalities above hold, and we have $|F(Q[S])|=3$. 
Thus $Q[S]$ has the unique $4$-face and two $6$-faces under our assumption, 
and hence there exists a $6$-face of $Q[S]$ that contains a connected component of 
$G-S$. 
If $|E(Q[S])|\leq 7$, then we have $|F(Q[S])|=2$ since $Q[S]$ has at least two faces with size at least $6$. 
In the same way, we can easily find our desired $6$-face; 
note that $|E(Q[S])|\geq 6$ by our assumption.   
\end{Proof}


Now we prove our second main result in the paper. 
\bigskip

\begin{Proofof}{Theorem \ref{thm:4,5-conn}}
First, we show (i) in the statement. 
By Theorem~\ref{thm:4-conn-gene}, 
the sufficiency is trivial , and hence we discuss 
the necessity below. 
Let $S$ be a $4$-cut of an O1TG $G$. 
By Lemma~\ref{lm:4cycle}, $Q[S]$ contains a $4$-gonal $2$-cell face 
containing a vertex since the surface is $T^2$ now; observe 
that $T^2$ is orientable, and that any surface separating cycle is 
trivial on $T^2$.  
Therefore, (i) of the statement holds. 

Next, we show (ii). 
First, we discuss the sufficiency. 
We put $C_1=xy_1z_1x$ and $C_2=xy_2z_2x$, which are two homotopic cycles in the statement. 
By Proposition~\ref{prop:parity}, both $C_1$ and $C_2$ are essential. 
First, consider the $2$-cell region $R$ bounded by $C_1\cup C_2$. 
We may assume that the boundary closed walk of $R$ is $xy_1z_1xz_2y_2x$. 
If $R$ contains no vertex of $G$, then $R$ has the unique 
diagonal edge of $Q(G)$, say $y_1z_2$ without loss of generality; 
that is, $R$ contains exactly two faces of $Q(G)$. 
However in this case, a crossing edge $xz_2$ in $R$ and a non-crossing edge 
$xz_2$ on the boundary of $R$ would form multiple edges, a contradiction. 
Similarly, the other region, which is outside of $R$, must contain a vertex 
of $G$. 
Therefore, $\{x, y_1, y_2, z_1, z_2\}$ is a $5$-cut of $G$. 

Secondly, 
we discuss the necessity, by putting a minimal $5$-cut 
$S=\{x, y_1, y_2, z_1, z_2\}$. 
By Lemma~\ref{lm:5conn}, 
$Q[S]$ has a $6$-face $f$ such that 
$V(\partial f)=S$ and $f$ contains a connected component of $G-S$. 
Similar to the proof of Lemma~\ref{lm:4cycle}, 
we consider the following two cases: 
Case (i) $\partial f$ has at least two boundary components, 
and Case (ii) $\partial f$ has exactly one boundary component. 
First,  
we discuss Case (i). 
In this case, we may assume that $\partial f_1=W_1\cup W_2$ with 
$V(W_1)=\{x, y_1, y_2\}$ and $V(W_2)=\{x, z_1, z_2\}$; 
note that $\partial f$ has exactly two connected components, 
otherwise $f$ would not become a $6$-face. 
Then, we clearly have our desired two $3$-cycles; by Proposition~\ref{prop:parity} and by our former arguments. 
Under the situation, observe that the region bounded by 
$xy_1y_2xyz_2z_1x$ other than $f$ is homeomorphic to a $2$-cell, 
where $W_1=xy_1y_2x$ and $W_2=xz_1z_2x$ are homotopic to each other 
in this direction. 

Next, we consider Case (ii). 
We may assume that only $x$ appears exactly twice on $\partial f$; 
since $|V(\partial f)|=5$ now. 
By Lemma~\ref{lm:degree} and the simplicity of $G$, $x$ must appear 
on the diagonal position of $\partial f$. 
Then, we may assume that $\partial f=xy_1y_2xyz_2z_1x$, and 
hence we have our desired two $3$-cycles also in this case; 
indeed using Proposition~\ref{prop:parity}.   
\end{Proofof}


\section{Matching extendability of O1TGs}\label{sect:Proofs}

Let $G$ be an optimal $1$-embedded graph on $F^2$, and let $W$ 
be a closed walk consisting of only non-crossing edges that bounds a $2$-cell region $D$. 
If $D$ contains an odd number of vertices, then 
we call $D$ an {\em odd weighted region\/}. 
In particular, if $W$ is a cycle, then $W$ is a {\em barrier cycle\/}. 
A barrier cycle of length $k$ is called a {\em barrier $k$-cycle\/}. 

The following lemma is a generalization of Lemma~2.3 in \cite{P1} that 
are often used in matching theory; this can be easily 
derived from Tutte's $1$-factor theorem. 

\begin{lm}[\cite{FSS}]\label{lem:blocker}
Let $G$ be a graph of even order and 
$\{e_1, \ldots , e_{m+1}\}$ 
be a matching of $G$ which is not extendable. 
Then there exists $S \subseteq V(G)$ 
such that 
\begin{itemize}
\item[{(i)}]
$\displaystyle S \supset \bigcup_{i=1}^{m+1} V(e_i)$ and
\item[{(ii)}]	$|S| \leq C_o(G - S) + 2m$. 
\end{itemize}
\end{lm}


In fact, by assuming $m$-extendability of a graph $G$, 
the equality of (ii) in Lemma~\ref{lem:blocker} holds (see e.g., \cite{FSS}). 
In the remaining part of the section, 
we prove the following three main results concern with 
matching extendability, using tools proved in 
the former sections. 

\bigskip

\begin{Proofof}{Theorem \ref{thm:1-ext}}
First, we show the necessity. 
Suppose that an O1TG $G$ contains a subgraph $H$ in the statement. 
Since $H$ is a quadrangulation of $T^2$, 
$|V(H)|=|F(H)|$ holds, and hence we have $|V(H)|=C_o(G-H)$ 
by the assumption. 
Let $e$ be an edge in $H$. 
Then $G'=G-V(e)$ has a cut 
$V(H)\setminus V(e)$ 
such that $G'-(V(H)\setminus V(e))$ has exactly 
$|V(H)|$ odd components. 
Thus $G'$ does not have a perfect matching by Tutte's $1$-factor theorem. 
That is, $G$ is not $1$-extendable. 

Next, we discuss the sufficiency. 
Let $G$ be an O1TG that is not $1$-extendable, and assume that $e$ is an edge of $G$ 
that is not extendable. 
Then there exists $S \subset V(G)$ such that 
(i) $V(e)\subset S$ and (ii) $|S| \leq C_o(G-S)$ 
by Lemma~\ref{lem:blocker} for $m=0$. 
Now we consider $Q[S]$ on $T^2$. 
By Lemma~\ref{lm:edge-bound}, we have $2|F(Q[S])| \geq |E(Q[S])|$ 
with $m=0$ and $\chi(F^2)=0$. 
On the other hand, 
$|E(Q[S])| \geq 2|F(Q[S])|$ holds by Lemma~\ref{lm:parity}, 
and hence we obtain $2|F(Q[S])|=|E(Q[S])|$. 
Actually, Lemma~\ref{lm:edge-bound} is obtained by using Euler's formula, 
and this equality implies that $Q[S]$ is a $2$-cell embedding, and 
it further indicates that $Q[S]$ is a quadragulation of $T^2$. 
It is well-known that the number of vertices equals the number 
of faces of a quadrangulation of $T^2$, and hence we have $|S|=|F(Q[S])|$. 
Then, $|F(Q[S])| \leq C_o(G-S)$ by (ii) above, 
and thus $|F(Q[S])|=C_o(G-S)$ by Lemma~\ref{lm:sep}. 
Hence each face of $Q[S]$ contains the unique odd component of $G-S$, 
that is, each face of $Q[S]$ is an odd weighted region. 
Therefore $Q[S]$ is a subgraph $H$ in the statement. 
\end{Proofof}



\bigskip

\begin{Proofof}{Theorem \ref{thm:2-ext}}
Any two independent edges on a barrier $4$-cycle are not extendable, 
and hence the necessity clearly holds. 
Thus, we discuss the sufficiency of the statement below. 
Let $G$ be an O1TG that is not $2$-extendable, and assume that 
$e_1$ and $e_2$ are independent edges of $G$ that are not extendable. 
Then there exists $S \subset V(G)$ such that  
(i) $V(e_1)\cup V(e_2) \subset S$ and (ii) $|S| \leq C_o(G-S)+2$ 
by Lemma~\ref{lem:blocker} for $m=1$. 

Now we consider $Q[S]$ on $T^2$. 
By Lemma~\ref{lm:edge-bound} with $m=1$, 
we have $2|F(Q[S])|+2 \geq |E(Q[S])|$. 
This inequality with Lemma~\ref{lm:parity} implies that there 
exist at most two faces of $Q[S]$ with size at least $6$; 
that is, the others are all $4$-faces. 
If $|S| \geq 5$, that is, $C_o(G-S) \geq 3$, 
then there exists a $4$-gonal face $f$ of $Q[S]$ 
that contains an odd component of $G-S$. 
Since any surface separating closed curve on $T^2$ is 
trivial, and hence, $\partial f$ is our desired barrier $4$-cycle. 
On the other hand, if $|S|=4$, then 
there is also a barrier $4$-cycle by (II) of Lemma~\ref{lm:4cycle} 
since $T^2$ is an orientable closed surface. 
Thus, we are done. 
\end{Proofof}


We use the following lemma in the proof of Theorem~\ref{thm:3-ext}. 

\begin{lm}\label{lem:oddcompo.}
Let $G$ be an $2q$-connected optimal $1$-embedded graph on $F^2$ with $q \geq3$ and 
$\{e_1, \ldots , e_{m+1}\}$ 
be a matching of $G$ which is not extendable. 
Then we have the following: 
$$C_o(G-S) \leq \displaystyle \frac{2m-\chi(F^2)}{q-2}$$
\end{lm}

\begin{Proof}
Then there exists $S \subset V(G)$ such that (i) 
$S \supset \bigcup_{i=1}^{m+1}V(e_i)$ and (ii) 
$|S| \leq C_o(G - S) + 2m$ of Lemma~\ref{lem:blocker}. 
Now we consider $Q[S]$ on $F^2$. 
Since $G$ is $2q$-connected, $Q[S]$ has at least 
$C_o(G-S)$ faces with size at least $2q$ 
by Lemma~\ref{lm:sep}. 
By (ii) of Lemma~\ref{lm:Q[S]}, 
$|S|-\chi(F^2) + (2- q)C_0(G-S) \geq |F(Q[S])|$ holds, 
and thus $C_o(G-S) +2m -\chi(F^2) + (2- q)C_0(G-S) 
= 2m -\chi(F^2) + (3- q)C_0(G-S) \geq |F(Q[S])|$. 
Since $|F(Q[S])| \geq C_o(G-S)$ by Lemma~\ref{lm:sep}, 
we have 
$2m -\chi(F^2) \geq  (q -2 )C_0(G-S)$, 
and hence the inequality in the statement holds. 
\end{Proof}


\begin{Proofof}{Theorem~\ref{thm:3-ext}}
The necessity is trivial and hence we prove the sufficiency of the statement; recall the observation in the introduction. 
Let $G$ be an $8$-regular O1TG. 
Suppose, for a contradiction, that $G$ is not $3$-extendable, and assume that 
$M=\{e_1, e_2, e_3\}$ is a matching that are not extendable. 
Then there exists $S \subset V(G)$ 
such that (i) $V(M) \subset S$ and (ii) $|S| \leq C_o(G - S) + 4$ by Lemma~\ref{lem:blocker} for $m=2$. 
Since $G$ is $8$-regular, $\kappa(G) \in \{6,8\}$ by Theorem~\ref{thm:8-conn}. 
First, suppose that $\kappa(G)=8$. 
Then the inequality 
$C_o(G-S) \leq 2$ holds by Lemma~\ref{lem:oddcompo.} with $q=4$, $m=2$ and 
$\chi(T^2)=0$, and hence 
we have $|S|\leq 6$ by (ii) above. This contradicts $\kappa(G)=8$. 

Hence we assume that $\kappa(G)=6$ and hence 
$Q(G)$ is isomorphic to 
$Q(3, r, q)$ $(q \geq 4$ and $r \geq 0)$ by 
Theorem~\ref{thm:8-conn}. 
Then the inequality 
$C_o(G-S) \leq 4$ holds by Lemma~\ref{lem:oddcompo.}, 
with $q=3$, $m=2$ and $\chi(T^2)=0$, and hence 
$|S|\leq 8$ holds by (ii) above.  
First, we consider the case of $|S|=6$. 
By Lemma~\ref{lm:6cut3}, $Q[S]$ is the union of two essential $3$-cycles $C_1$ and $C_2$ 
that are homotopic to each other. 
Since $S=V(M)$, one edge in $M$ 
joins vertices of $C_1$ and $C_2$. 
This implies that one of the two annular faces of $Q[S]$ does not contain a vertex of $G$, 
contradicting Lemma~\ref{lm:sep}. 

Next, we discuss the case of $|S|=7$.  
By Lemma~\ref{lm:min7-cut}, 
$S$ contains a $6$-cut $S'$, and we put $S\setminus S' = \{v\}$.  
By Lemma~\ref{lm:6cut3} again, $Q[S']$ is the union of two essential $3$-cycles. 
Let $H_1$ and $H_2$ be connected components of $G-S'$. 
Under the situation, 
$Q[V(H_i)]$ is a Cartesian product of a $3$-cycle and a $k$-path $(k\geq 0)$ for 
each $i \in \{1, 2\}$. 
Then each of $H_1'$ and $H_2'$ is either a $3$-cycle (if $k=0$) or $3$-connected (if $k\geq 1$). 
Thus, $H_i - v$ is connected for each $i \in \{1, 2\}$. 
Therefore $C_o(G-S)\leq 2$ holds, contradicting $C_o(G-S)\geq 3$ obtained from (ii) above. 
Finally, we consider the case of $|S|=8$. 
Similarly, by (ii) above, we have $C_o(G-S)\geq 4$, that is, $C_o(G-S)= 4$. 
By the argument in the proof of Lemma~\ref{lem:oddcompo.}, 
all the equalities hold, and hence 
$Q[S]$ is a $2$-cell embedding, and has exactly four faces, 
each of which is $6$-gonal, and contains an odd component of $G-S$ . 
This clearly contradicts (i) of Lemma~\ref{lm:4regu.}.  
\end{Proofof}


\section{Remarks}\label{sect:remarks}
In this paper, we have discussed connectivity and matching extendability of 
O1TGs. 
We aim to investigate similar problems on optimal $1$-embedded graphs on the Klein bottle, 
and further extend our discussion to such graphs on 
general closed surfaces. 
Actually, we proved some results that works for optimal $1$-embedded graphs 
on general closed surfaces, for example, Theorem~\ref{thm:4-conn-gene}, 
Lemmas~\ref{lm:4cycle} and \ref{lem:oddcompo.}. 
Furthermore, there are facts that do not hold for general closed surfaces but do apply to 
optimal $1$-embedded graphs on the torus and the Klein bottle, 
that is, closed surfaces with Euler characteristic $0$.  
For example, Lemmas~\ref{lm:6-conn.}, \ref{lm:6cut2} and \ref{lm:5conn} are among them, 
and further, we can confirm that the assertion of Theorem~\ref{thm:1-ext} also 
holds for optimal $1$-embedded graphs on the Klein bottle. 
These results are expected to be very helpful for our future research. 
However, on the other hand, there are many claims that require us to re-examine 
the discussion using properties specific to the Klein bottle. 
Specifically, it is known that the Klein bottle has 
a non-trivial surface separating simple closed curve, 
which is known as an {\em equator\/}, 
and understanding their impact on our problem seems to be the first step of our 
future works. 

We believe that, at least, optimal $1$-embedded graphs on the Klein bottle has no 
minimal $7$-cut. 
In other words, we propose the following conjecture, 
which is similar to the proposition for O1TGs. 

\begin{conj}\label{conj1}
Every $8$-regular optimal $1$-embedded graph $G$ on the Klein bottle has connectivity either $6$ or $8$. 
\end{conj}

Similar to the case of the torus, $4$-regular quadrangulations on the Klein bottle 
are well classified in \cite{NakaNega}, 
and it would be even better if we could describe our claims using those classifications.
In the end of the paper, we furhter establish the following conjecture 
concern with matching extendability of those graphs.  

\begin{conj}\label{conj2}
Every $8$-regular optimal $1$-embedded graph on the Klein bottle is 
$3$-extendable. 
\end{conj}

\bigskip

\noindent
{\bf Statements and Declarations}
\bigskip

This work was supported by JSPS KAKENHI Grant Number 23K03196 and JST SPRING, Grant Number JPMJSP2121. 
The authors have no relevant financial or non-financial interests to disclose.

%



\end{document}